\documentclass[12pt]{article}
\usepackage{latexsym}
\usepackage{amsmath}
\usepackage{amssymb}

\def \zta {\mbox{\boldmath $\zeta$}}

\begin{document}

\thispagestyle{empty}

\centerline{\Large{\bf  Extreme Points of the Convex Set of Joint Probability}}
\centerline{\Large{\bf  Distributions with Fixed Marginals}}
\centerline{\bf by}
\centerline{\bf K. R. Parthasarathy}
\centerline{\bf Indian Statistical Institute, Delhi Centre,}
\centerline{\bf 7, S. J. S. Sansanwal Marg, }
\centerline{\bf New Delhi - 110 016, India.}
\centerline{\bf e-mail : krp@isid.ac.in}
\vskip40pt

\vskip 0.5in
\noindent {\bf Summary :}  By using a quantum probabilistic approach we obtain a description of the extreme points of the convex set of all joint probability distributions on the product of two standard Borel spaces with fixed marginal distributions.

\vskip 0.5in
\noindent AMS   Subject Classification Index 46L53, 15A51

\newpage
\section{Introduction}
It is a well-known theorem of Garret Birkhoff \cite{gb} and von Neumann \cite{jvn}, \cite{rbbtesr}, \cite{rbh} that the extreme points in the convex set of all $n \times n$ bistochastic (or doubly stochastic) matrices are precisely the $n$-th order permutation matrices. Here we address the following problem: If $G$ is a standard Borel group acting measurably on two standard probability spaces $(X_i, {\mathcal F}_i, \mu_i),$ $i=1,2$ where $\mu_i$ is invariant under the $G$-action for each $i$ then what are the extreme points of the convex set of all joint probability distributions on the product Borel space $(X_1 \times X_2, {\mathcal F}_1 \otimes {\mathcal F}_2)$ which are invariant under the diagonal action $(x_1, x_2) \mapsto (gx_1, gx_2)$ where $x_i \in X_i,$ $i=1,2$ and $g \in G?$

Our approach to the problem mentioned above is based on a quantum probabilistic method arising from Stinespring's \cite{wfs} description of completely positive maps on $C^{\ast}$ algebras. We obtain a necessary and sufficient condition for the extremality of a joint distribution in the form of a regression condition. This leads to examples of extremal nongraphic joint distributions in the unit square with uniform marginal distributions on the unit interval. The Birkhoff-von Neumann theorem is deduced as a corollary of the main theorem.

\section{The convex set of covariant bistochastic maps on $C^{\ast}$ algebras}

For any complex separable Hilbert space ${\mathcal H},$ express its scalar product in the Dirac  notation $\langle \cdot | \cdot \rangle $ and denote by ${\mathcal B}({\mathcal H})$ the $C^{\ast}$ algebra of all bounded operators on ${\mathcal H}.$ Let $G$ be a group with fixed unitary representations $g \mapsto U_g,$ $g \mapsto V_g,$ $g \in G$ in Hilbert spaces ${\mathcal H}_1,$ ${\mathcal H}_2$ respectively and let ${\mathcal A}_i \subset {\mathcal B} ({\mathcal H}_i),$ $i = 1,2$ be unital $C^{\ast}$ algebras invariant under respective conjugations by $U_g,$ $V_g$ for every $g$ in $G.$ Let  $\omega_i$ be a fixed state in ${\mathcal A}_i$ for each $i,$ satisfying the invariance conditions:
\begin{equation}
\omega_1 \left ( U_g X U_g^{-1}\right ) = \omega_1 (X), \omega_2 (V_g Y V_g^{-1}) =\omega_2 (Y) \quad \forall \,\,\, X \in {\mathcal A}_1, \,\,Y \in {\mathcal A}_2, g \in G. \label{eq2.1}
\end{equation}
Consider a linear, unital and completely positive map $T: {\mathcal A}_1 \rightarrow {\mathcal A}_2$ satisfying the following:
\begin{eqnarray}
 \omega_2 (T(X)) &=& \omega_1 (X) \quad \forall \,\,\, X \in {\mathcal A}_1, \label{eq2.2}\\
T \left (U_g X U_g^{-1} \right ) &=& V_g T(X) V_g^{-1} \quad \forall \,\,\,X \in {\mathcal A}_1, g \in G. \label{eq2.3}
\end{eqnarray}
Then we say that $T$ is a $G$-{\it covariant bistochastic map} with respect to the pair of states $\omega_1, \omega_2$ and representations $U.,$ $V.$. Denote by $\mathbb{K}$ the convex set of all such covariant bistochastic maps from ${\mathcal A}_1$ into ${\mathcal A}_2.$ We shall now present a necessary and sufficient condition for an element $T$ in $\mathbb{K}$ to be an extreme point of $\mathbb{K}.$

To any $T \in \mathbb{K}$ we can associate a Stinespring triple $({\mathcal K}, j, \Gamma)$ where ${\mathcal K}$ is a Hilbert space, $j$ is a $C^{\ast}$ homomorphism from ${\mathcal A}_1$ into ${\mathcal B}({\mathcal K})$ and $\Gamma$ is an isometry from ${\mathcal H}_2$ into ${\mathcal K}$ satisfying the following properties:
\begin{itemize}
 \item [(i)] $\Gamma^{\dagger} j (X) \Gamma = T(X) \quad \forall \,\,X \in {\mathcal A}_1;$
\item [(ii)] The linear manifold generated by $\left \{j (X) \Gamma u \big | u \in {\mathcal H}_2, X \in {\mathcal A}_1 \right \}$ is dense in ${\mathcal K}.$
\end{itemize}
Such a Stinespring triple is unique upto a unitary isomorphism, i.e., if $({\mathcal K}^{\prime},j^{\prime}, \Gamma^{\prime})$ is another triple satisfying the properties (i) and (ii) above then there exists a unitary isomorphism $\theta : {\mathcal K} \rightarrow {\mathcal K}^{\prime}$ such that $\theta j (X) = j^{\prime} (X) \theta$ $\forall$ $X \in {\mathcal A}_1$ and $\theta \Gamma v = \Gamma^{\prime} v$ $\forall$ $v \in {\mathcal H}_2.$ (See \cite{wfs}.)

We now claim that the covariance property of $T$ ensures the existence of a unitary representation $g \mapsto W_g$ of $G$ in ${\mathcal K}$ satisfying the relations:
\begin{eqnarray}
W_g j (X) \Gamma u &=& j (U_g X U_g^{-1})\Gamma V_g u \quad \forall \,\, X \in {\mathcal A}_1, g \in G, u \in {\mathcal H}_2, \, \label{eq2.4} \\
W_g j (X) W_g^{-1} &=& j (U_g X U_g^{-1}) \quad \forall \,\,X \in {\mathcal A}_1, g \in G. \label{eq2.5}
\end{eqnarray}
Indeed, for any $X,Y$ in ${\mathcal A}_1$ $u,$ $v \in {\mathcal H}_2$ and $g \in G$ we have from the properties (i) and (ii) above and \eqref{eq2.3}
\begin{eqnarray*}
\lefteqn{ \langle  j \left ( U_g X U_g^{-1} \right ) \Gamma V_g u \big | j   \left (U_g Y U_g^{-1} \right ) \Gamma V_g v \rangle } \\
&=& \langle u \big | V_g^{-1} \Gamma^{\dagger} j \left (U_g X^{\dagger} Y U_g^{-1} \right ) \Gamma V_g v \rangle          \\
&=&  \langle u \big | V_g^{-1} T (U_g X^{\dagger} Y U_g^{-1}) V_g v \rangle \\ 
&=& \langle u \big | T (X^{\dagger} Y) \big | v \rangle \\
&=& \langle j (X) \Gamma u \big | j (Y) \Gamma v \rangle.
\end{eqnarray*}
In other words, the correspondence $j (X) \Gamma u \mapsto j (U_g X U_g^{-1}) \Gamma V_g u$ is a scalar product preserving map on a total subset of ${\mathcal K},$ proving the claim.
\vskip0.2in
\noindent{\bf Theorem 2.1 \quad}  Let $T \in \mathbb{K}$ and let $({\mathcal K},j,\Gamma)$ be a Stinespring triple associated to $T.$  Let $g \mapsto W_g$ be the unique unitary representation of $G$ satisfying the relations \eqref{eq2.4} and \eqref{eq2.5}. Then $T$ is an extreme point of $\mathbb{K}$ if and only if there exists no nonzero hermitian operator $Z$ in the commutant of the set $\left \{j (X), X \in {\mathcal A}_1 \right \} \cup \left \{ W_g, g \in G \right \}$ satisfying the following two conditions:
\begin{itemize}
 \item [(i)] $\Gamma^{\dagger} Z \Gamma = 0;$
\item [(ii)] $\Gamma^{\dagger} Z j (X) \Gamma  \in {\mathcal A}_2$ and $\omega_2 \left (\Gamma^{\dagger} Z j (X) \Gamma \right ) = 0$ $\qquad \forall$ $X \in {\mathcal A}_1.$
\end{itemize}
\vskip0.2in
\noindent{\bf Proof \quad} Suppose $T$ is not an extreme point of $\mathbb{K}.$ Then there exist $T_1, T_2 \in \mathbb{K},$ $T_1 \neq T_2$ such that $T = \frac{1}{2} (T_1 + T_2).$ Let $({\mathcal K}_1, j_1, \Gamma_1)$ be a Stinespring triple associated to $T_1.$ Then by the argument outlined in the proof of Proposition 2.1 in \cite{krp} there exists a bounded operator $J : {\mathcal K} \rightarrow {\mathcal K}_1$ satisfying the following properties:
\begin{itemize}
 \item [(i)] $J j (X) \Gamma u = j_1 (X) \Gamma_1 u$ $\qquad \forall$ $X \in {\mathcal A}_1,$ $u \in {\mathcal H}_2;$
\item [(ii)] The positive operator $\rho:= J^{\dagger}J$ is in the commutant of $\left \{ j (X), X \in {\mathcal A} \right \}$ in ${\mathcal B}({\mathcal K});$
\item [(iii)] $T_1 (X) = \Gamma^{\dagger} \rho j (X) \Gamma.$
\end{itemize}
Since $T_1 \neq T_2$ it follows that $T_1 \neq T$ and hence $\rho$ is different from the identity operator. We now claim that $\rho$ commutes with $W_g$ for every $g$ in $G.$ Indeed, for any $X,Y$ in ${\mathcal A}_1,$ $u, v$ in ${\mathcal H}_2$ we have from the definition of $\rho$ and $J,$ equation \eqref{eq2.4} and the covariance of $T_1$
\begin{eqnarray*}
\lefteqn{ \langle  j(X) \Gamma u \big |\rho W_g   \big | j(Y) \Gamma v  \rangle } \\
&=& \langle  j(X) \Gamma u  \big |J^{\dagger} J \big | j (U_g Y U_g^{-1}) \Gamma V_g v \rangle\\
&=& \langle  j_1 (X) \Gamma_1 u \big |j_1 (U_g YU_g^{-1}) \Gamma_1 V_g v   \rangle   \\
&=& \langle  u \big | \Gamma_1^{\dagger} j_1 (X^{\dagger} U_g Y U_g^{-1})\Gamma_1  \big |V_g v  \rangle   \\
&=& \langle  u \big | T_1 (X^{\dagger} U_g Y U_g^{-1})  \big |V_g v   \rangle   \\
&=& \langle  u \big | V_g T_1 (U_g^{-1} X^{\dagger} U_g Y)   \big |v  \rangle.
\end{eqnarray*}
On the other hand, by the same arguments, we have
\begin{eqnarray*}
\lefteqn{ \langle j (X) \Gamma u   \big |W_g \rho   \big | j (Y) \Gamma v  \rangle} \\
&=& \langle  j (U_g^{-1} X U_g) \Gamma V_g^{-1} u \big |J^{\dagger} J \big | j (Y) \Gamma v  \rangle   \\
&=& \langle j_1 (U_g^{-1} X U) \Gamma_1 V_g^{-1} u   \big |j_1 (Y) \Gamma_1 v     \rangle   \\
&=& \langle   u \big | V_g T_1 (U_g^{-1} X^{\dagger} U_g Y)  \big |v  \rangle   
\end{eqnarray*}
Comparing the last two identities and using property (ii) of the Stinespring triple we conclude that $\rho$ commutes with $W_g.$ Putting $Z = \rho - I$ we have
\begin{equation}
\Gamma^{\dagger} Z j (X) \Gamma = T_1 (X) - T(X) \quad \forall \,\,\, X \in {\mathcal A}_1. \label{eq2.6}
\end{equation}
Clearly, the right hand side of this equation is an element of ${\mathcal A}_2$ and
$$\omega_2 (\Gamma^{\dagger} Z j (X) \Gamma) = \omega_1 (X) - \omega_1 (X) = 0 \quad \forall \,\,\, X \in {\mathcal A}_1. $$
Putting $X=I$ in \eqref{eq2.6} we have $\Gamma^{\dagger} Z \Gamma = 0.$ Then $Z$ satisfies properties (i) and (ii) in the statement of the theorem, proving the sufficiency part.

Conversely, suppose there exists a nonzero hermitian operator $Z$ in the commutant of $\left \{ j (X), X \in {\mathcal A}_1 \right \} \cup \left \{W_g, g \in G \right \}$ satisfying properties (i) and (ii) in the theorem. Choose and fix a positive constant $\varepsilon$ such that the operators $I \pm \varepsilon Z$ are positive. Define the maps $T_{\pm} : {\mathcal A}_1 \rightarrow {\mathcal A}_2$ by
\begin{equation}
T_{\pm} (X) = \Gamma^{\dagger} (I \pm \varepsilon Z) j (X) \Gamma, \qquad  X \in {\mathcal A}_1. \label{eq2.7}
\end{equation}
Since
$$(I \pm \varepsilon Z)j(X) = \sqrt{I \pm \varepsilon Z} j (X) \sqrt{I \pm \varepsilon Z} $$
it follows that $T_{\pm}$ are completely positive. By putting $X = I$ in \eqref{eq2.7} and using property (i) of $Z$ in the theorem we see that $T_{\pm}$ are unital. Furthermore, we have from equations \eqref{eq2.4} and \eqref{eq2.5}, for any $g \in G,$ $X \in {\mathcal A}_1,$
\begin{eqnarray*}
 T_{\pm} (U_g X U_g^{-1}) &=& \Gamma^{\dagger} (I \pm \varepsilon Z) W_g j (X) W_g^{-1} \Gamma \\
&=& V_g \Gamma^{\dagger} (I \pm \varepsilon Z) j (X) \Gamma V_g^{-1} \\
&=& V_g T_{\pm} (X) V_g^{-1}.
\end{eqnarray*}
Also, by property (ii) in the theorem we have
$$\omega_2 (T_{\pm}(X)) = \omega_2 (T(X)) = \omega_1 (X) \quad \forall \,\,\, X \in {\mathcal A}_1. $$
Thus $T_{\pm} \in \mathbb{K}.$ Note that
$$\langle u \big | \Gamma^{\dagger} Z j (X^{\dagger}Y) \Gamma \big | v  \rangle = \langle j (X) \Gamma u \big | Z \big | j (Y) \Gamma v   \rangle $$
cannot be identically zero when $X$ and $Y$ vary in ${\mathcal A}_1$ and $u$ and $v$ vary in ${\mathcal H}_2.$ Thus $\Gamma^{\dagger} Z j (X) \Gamma \not\equiv 0$ and hence  $T_{+} \neq T_{-}.$ But $T = \frac{1}{2} (T_{+} + T_{-}).$  In other words $T$ is not an extreme point of $\mathbb{K}.$ This proves necessity. $\square$

\section{The convex set of invariant joint distributions with fixed marginal distributions}
\setcounter{equation}{0}

Let $(X_i, {\mathcal F}_i, \mu_i),$ $i = 1,2$ be standard probability spaces and let $G$ be a standard Borel group acting measurably on both $X_1$ and $X_2$ preserving $\mu_1$ and $\mu_2.$ Denote by $\mathbb{K} (\mu_1, \mu_2)$ the convex set of all joint probability distributions on the product Borel space $(X_1, \times X_2,$ \\$ {\mathcal F}_1 \otimes {\mathcal F}_2)$ invariant under the diagonal $G$ action $(g, (x_1, x_2)) \mapsto (gx_1, gx_2),$ $x_i \in X_i,$ $g \in G$ and having the marginal distribution $\mu_i$ in $X_i$ for each $i.$ Choose and fix $\omega \in \mathbb{K} (\mu_1, \mu_2).$ Our present aim is to derive from the quantum probabilistic result in Theorem 2.1, a necessary and sufficient condition for $\omega$ to be an extreme point of $\mathbb{K}(\mu_1, \mu_2).$ To this end we introduce the Hilbert spaces ${\mathcal H}_i = L^2 (\mu_i),$ ${\mathcal K} = L^2 (\omega)$ and the abelian von Neumann algebras ${\mathcal A}_i \subset {\mathcal B}({\mathcal H}_i)$ where ${\mathcal A}_i = L^{\infty} (\mu_i)$ is also viewed as the algebra of operators of multiplication by functions from $L^{\infty}(\mu_i).$ For any $\varphi \in L^{\infty} (\mu_i)$ we shall denote by the same symbol  $\varphi$ the multiplication operator $f \mapsto \varphi f,$ $f \in L^2 (\mu_i).$ For any $\varphi \in {\mathcal A}_1$ define the operator $j (\varphi)$ in ${\mathcal K}$ by
\begin{equation}
(j (\varphi)f) (x_1, x_2) = \varphi (x_1) f (x_1, x_2), \quad f \in {\mathcal K}, x_i \in X_i. \label{eq3.1}
\end{equation}
Then the correspondence $\varphi \mapsto j (\varphi)$ is a von Neumann algebra homomorphism from ${\mathcal A}_1$ into ${\mathcal B}({\mathcal K}).$ Define the isometry $\Gamma : {\mathcal H}_2 \rightarrow {\mathcal K}$ by
\begin{equation}
(\Gamma v) (x_1, x_2) = v (x_2), \quad v \in {\mathcal H}_2. \label{eq3.2}
\end{equation}
Then, for $f \in {\mathcal K},$ $v \in {\mathcal H}_2$ we have
\begin{eqnarray*}
\langle f \big | \Gamma v \rangle &=& \int_{X_1 \times X_2} \bar{f} (x_1, x_2) v (x_2) \omega (dx_1 dx_2) \\
&=& \int_{X_{2}} \mu_2 (dx_2) \left [\bar{f} (x_1, x_2) \nu (dx_1, x_2)  \right ] v (x_2)
\end{eqnarray*}
where $\nu (E, x_2),$ $E \in {\mathcal F}_1,$ $x_2 \in X_2$ is a measurable version of the conditional probability distribution on ${\mathcal F}_1$ given the sub $\sigma$-algbera $\{X_1 \times F, F \in {\mathcal F}_2 \} \subset {\mathcal F}_1 \otimes {\mathcal F}_2.$ Thus the adjoint $\Gamma^{\dagger} : {\mathcal K} \rightarrow {\mathcal H}_2$ of $\Gamma$ is given by
\begin{equation}
(\Gamma^{\dagger} f) (x_2) = \int_{X_{1}} f (x_1, x_2) \nu (dx_1, x_2). \label{eq3.3}
\end{equation}
Hence
\begin{eqnarray}
 (j (\varphi) \Gamma v)(x_1, x_2) &=& \varphi (x_1) v (x_2), \qquad \varphi \in {\mathcal A}_1, \,\,\, v \in {\mathcal H}_2, \label{eq3.4} \\
(\Gamma^{\dagger} j (\varphi) \Gamma v)(x_2) &=& \left [ \int \varphi (x_1) \nu (dx_1, x_2) \right ] v (x_2).          \label{eq3.5}
\end{eqnarray}
In other words
\begin{equation}
\Gamma^{\dagger} j (\varphi) \Gamma = T (\varphi) \label{eq3.6}
\end{equation}
where $T (\varphi) \in {\mathcal A}_2$ is given by
\begin{equation}
T(\varphi) (x_2) = \int_{X_{1}} \varphi (x_1) \nu (dx_1, x_2). \label{eq3.7}
\end{equation}
Equations \eqref{eq3.1}-\eqref{eq3.7} imply that $T$ is a linear, unital and positive (and hence completely positive) map from the abelian von Neumann algebra ${\mathcal A}_1$ into ${\mathcal A}_2$ and $({\mathcal K}, j, \Gamma)$ is, indeed, a Stinespring triple for $T.$ Furthermore, the unitary operators $U_g, V_g$ and $W_g$ in ${\mathcal H}_1, \,{\mathcal H}_2$ and ${\mathcal K}$ respectively defined by
\begin{eqnarray*}
(U_g u) (x_1) &=& u (g^{-1} x_1), \quad u \in {\mathcal H}_1,\\
(V_g v) (x_2) &=& v (g^{-1} x_2), \quad v \in {\mathcal H}_2, \\
(W_g f) (x_1, x_2) &=& f (g^{-1} x_1, \,g^{-1} x_2), \quad f \in k
\end{eqnarray*}
satisfy the relations \eqref{eq2.4} and \eqref{eq2.5}.

Our next lemma describes operators of the form $Z$ occurring in Theorem 2.1.
\vskip0.2in
\noindent{\bf Lemma 3.1 \quad} Let $Z$ be a bounded hermitian operator in ${\mathcal K}$ satisfying the following conditions:
\begin{itemize}
 \item [(i)] $Z j (\varphi) = j (\varphi) Z \quad \forall \,\,\, \varphi \in {\mathcal A}_1,$
\item [(ii)] $ZW_g = W_g Z \quad \forall \,\,\, g \in G,$
\item [(ii)] $\Gamma^{\dagger} Z j (\varphi) \Gamma \in {\mathcal A}_2 \quad \forall \,\,\, \varphi \in {\mathcal A}_1.$
\end{itemize}
Then there exists a function $\zeta \in L^{\infty} (\omega)$ satisfying the following properties:
\begin{itemize}
 \item [(a)] $\zeta (gx_1, gx_2 ) = \zeta (x_1, x_2) \,\mbox{a.e.}\, (\omega) \quad \forall \,\,\, g \in G,$
\item [(a)] $(Zf) (x_1, x_2) = \zeta (x_1, x_2) f (x_1, x_2) \qquad \forall \,\,\, f \in {\mathcal K}$
\end{itemize}
\vskip0.2in
\noindent{\bf Proof \quad} Let
$$\zeta (x_1, x_2) = (Z1)(x_1, x_2)$$
where the symbol $1$ also denotes the function identically equal to unity. For functions $u,v$ on $X_1, X_2$ respectively denote by $u \otimes v$ the function on $X_1 \times X_2$ defined by $u \otimes v (x_1, x_2) = u (x_1) v (x_2).$ By property (i) of $Z$ in the lemma we have
\begin{eqnarray}
(Z \varphi \otimes 1) (x_1, x_2) &=& (Z j (\phi)1) (x_1, x_2) \nonumber \\
&=& (j (\phi) Z1) (x_1, x_2)  \nonumber\\
&=& \varphi (x_1) \zeta (x_1, x_2) \qquad \forall \,\,\,\varphi \in {\mathcal A}_1. \label{eq3.8}
\end{eqnarray}
If $\varphi \in {\mathcal A}_1,$ $v \in {\mathcal H}_2,$ we have
\begin{eqnarray}
(Z \varphi \otimes v) (x_1, x_2) &=& (Z j (\varphi) \Gamma v) (x_1, x_2) \nonumber \\
&=& (j (\varphi) Z \Gamma v) (x_1, x_2) \nonumber \\
&=& \varphi (x_1) (Z 1 \otimes v) (x_1, x_2) \label{eq3.9} 
\end{eqnarray}
From properties (i) and (iii) of $Z$ in the lemma and equations \eqref{eq3.3}, \eqref{eq3.8} and \eqref{eq3.9} we have
\begin{eqnarray*}
(\Gamma^{\dagger} Z j (\varphi) \Gamma v)(x_2) &=& \int (Z \varphi \otimes v)\nu (dx_1, x_2) \\
&=& \int \varphi (x_1) (Z 1 \otimes v) (x_1,x_2) \nu (dx_1, x_2)
\end{eqnarray*}
whereas the left hand side is of the form $R(\varphi) (x_1) v (x_2)$ for some $R(\varphi) \in L^{\infty} (\mu_2).$ Thus
$$R (\varphi) (x_2) v (x_2) = \int \varphi (x_1) (Z 1 \otimes v) (x_1, x_2) \nu (dx_1, x_2).$$
Choosing $v =1 $ we have from the definition of $\zeta$
$$R (\varphi) (x_2) = \int \varphi (x_1) \zeta (x_1, x_2) \nu (dx_1, x_2).$$
Thus, for every $\varphi \in {\mathcal A}_1$
$$\int \varphi (x_1) \zeta (x_1, x_2) v (x_2) \nu (dx_1, x_2) = \int \varphi (x_1) (Z 1 \otimes v) (x_1, x_2) \nu (dx_1, x_2) $$
and hence
$$(Z 1 \otimes v) (x_1, x_2) = \zeta (x_1, x_2) v (x_2) \,\mbox{a.e.}\, x_1 (\nu (., x_2)) \,\mbox{a.e.}\, x_2 (\mu_2). $$
Applying $j (\varphi)$ on both sides we get
$$(Z \varphi \otimes v) (x_1, x_2) = \zeta (x_1, x_2) \varphi (x_1) v (x_2) \,\mbox{a.e.}\, (\omega). $$
In other words $Z$ is the operator of multiplication by $\zeta$ and it follows that $\zeta \in L^{\infty} (\omega).$ Now property (ii) of $Z$ implies property (a) in the lemma. $\square$
\vskip0.2in
\noindent{\bf Theorem 3.2 \quad} Let $\omega \in \mathbb{K} (\mu_1, \mu_2).$ Then $\omega$ is an extreme point of $\mathbb{K} (\mu_1, \mu_2)$ if and only if there exists no nonzero real-valued function $\zeta \in L^{\infty} (\omega)$ satisfying the following conditions:
\begin{itemize}
 \item [(i)] $\zeta (gx_1, gx_2) = \zeta (x_1, x_2)$ a.e. $\omega$ $\forall$ $g \in G;$
\item [(ii)] $\mathbb{E} (\zeta (\xi_1, \xi_2) \big | \xi_1) = 0,$ $\mathbb{E} (\zeta (\xi_1, \xi_2) \big | \xi_2) = 0$ where $(\xi_1, \xi_2)$ is an $X_1 \times X_2$-valued random variable with distribution $\omega.$
\end{itemize}
\vskip0.2in
\noindent{\bf Proof \quad} Let $Z$ be a bounded selfadjoint operator in the commutant of $\left \{j (\varphi), \varphi \in {\mathcal A}_1 \right \} \cup \left \{W_g, g \in G \right \}$ such that $\Gamma^{\dagger} Z j (\varphi) \Gamma \in {\mathcal A}_2$ $\forall$ $\varphi \in {\mathcal A}_1.$ Then by Lemma 3.1 it follows that $Z$ is of the form
$$(Zf) (x_1, x_2) = \zeta (x_1, x_2)f (x_1, x_2)$$
where $\zeta \in L^{\infty}(\omega)$ and $\zeta (gx_1, gx_2) = \zeta (x_1, x_2)$ a.e. $(\omega).$ Note that
$$(\Gamma^{\dagger} Z \Gamma v)(x_2) = \left [\int_{X_{1}} \zeta (x_1, x_2) \nu (dx_1, x_2) \right ] v (x_2) \,\,\mbox{a.e.}\,\, (\mu_2), v \in {\mathcal H}_2.$$
Thus $\Gamma^{\dagger} Z \Gamma = 0$ if and only if $\mathbb{E} (\zeta (\xi_1, \xi_2) \big | \xi_2) = 0.$ Now we evaluate
$$(\Gamma^{\dagger} Z j (\varphi)\Gamma v )(x_2) = \int \varphi (x_1) v (x_2) \zeta (x_1,x_2) \nu (dx_1, x_2) \quad \mbox{a.e.}\,\, (\mu_2).$$
Looking upon $\Gamma^{\dagger} Z j (\varphi) \Gamma$ as an element of ${\mathcal A}_2$ and evaluating the state $\mu_2$ on this element we get
\begin{eqnarray*}
\mu_2 (\Gamma^{\dagger} Z j (\varphi) \Gamma) &=& \int \varphi (x_1) \zeta (x_1, x_2) \nu (dx_1, x_2) \mu (dx_2) \\
&=&  \int \varphi (x_1) \zeta (x_1, x_2) \omega (dx_1\, dx_2)           \\
&=& \mathbb{E}_{\omega} \varphi (\xi_1) \zeta (\xi_1, \xi_2) \\
&=& \mathbb{E}_{\mu_{1}} \varphi (\xi_1) \mathbb{E} (\zeta (\xi_1, \xi_2) \big | \xi_1).
\end{eqnarray*}
Thus $\mu_2 (\Gamma^{\dagger} Z j (\varphi) \Gamma) = 0$ $\forall$ $\varphi \in {\mathcal A}_1$ if and only if $\mathbb{E} (\zeta (\xi_1, \xi_2) \big | \xi_1) = 0.$ Now an application of Theorem 2.1 completes the proof of the theorem. $\square$

We shall now look at the special case when $G$ is the trivial group consisting of only the identity element. Let $(X_i, {\mathcal F}_i, \mu_i),$ $i=1,2$ be standard probability spaces and let $T : X_1 \rightarrow X_2$ be a Borel map such that $\mu_2 = \mu_1 T^{-1}.$ Consider an $X_1$-valued random variable $\xi$ with distribution $\mu_1.$ Then the joint distribution $\omega$ of the pair $(\xi, T \circ \xi)$ is an element of $\mathbb{K} (\mu_1, \mu_2)$ and by Theorem 2.1 is an extreme point. Similarly, if $T : X_2 \rightarrow X_1$ is a Borel map such that $\mu_2 T^{-1} = \mu_1$ and $\eta$ is an $X_2$-valued random variable with distribution $\mu_2$ then $(T \circ \eta, \eta)$ has a joint distribution which is an extreme point of $\mathbb{K} (\mu_1, \mu_2).$ Such extreme points are called {\it graphic} extreme points. Thus there arises the natural question whether there exist nongraphic extreme points. Our next lemma facilitates the construction of nongraphic extreme points.

\vskip0.2in
\noindent{\bf Lemma 3.3 \quad} Let $(X, {\mathcal F}, \lambda),$ $(Y, {\mathcal G}, \mu),$ $(Z, {\mathcal K}, \nu)$ be standard probability spaces and let $\xi, \eta, \zeta$ be random variables on a probability space with values in $X,Y,Z$ and distribution  $\lambda, \mu, \nu$ respectively. Suppose $\zeta$ is independent of $(\xi, \eta)$ and the joint distribution $\omega$ of $(\xi, \eta)$ is an extreme point of  $\mathbb{K} (\lambda, \mu).$ Let $\widetilde{\lambda},\widetilde{\mu}, \widetilde{\omega}$ be the distributions of $(\xi, \zeta),$ $(\eta, \zeta)$ and $((\xi, \zeta), (\eta, \zeta))$ respectively in the spaces $X \times Z,$  $Y \times Z$ and $(X \times Z) \times (Y \times Z).$ Then $\widetilde{\omega}$ is an extreme point of $\mathbb{K} (\widetilde{\lambda}, \widetilde{\mu}).$
\vskip0.2in
\noindent{\bf Proof \quad} Let $f$ be a bounded real-valued measurable function on $(X \times Z) \times (Y \times Z)$ satisfying the relations
\begin{eqnarray*}
\mathbb{E} \left \{f ((\xi, \zeta), (\eta, \zeta)) \big | (\eta, \zeta) \right \} &=& 0,\\
\mathbb{E} \left \{f ((\xi, \zeta), (\eta, \zeta)) \big | (\xi, \zeta) \right \} &=& 0.
\end{eqnarray*}
If we write
$$F_z (x,y) = f ((x,z), (y,z)) \quad \mbox{where} \quad (x,y,z) \in X \times Y \times Z$$
then we have
$$\mathbb{E} (F_z (\xi, \eta) \big | \eta) = 0, \quad \mathbb{E} (F_z (\xi, \eta) \big | \xi) = 0 \,\,\, \mbox{a.e.}\,\,z (\nu).  $$
Since $\omega$ is extremal it follows that $F_z (\xi, \eta) = 0$ a.e. $z (\nu)$ and therefore $f ((\xi,\zeta), (\eta, \zeta))=0.$ By Theorem 3.1 it follows that $\widetilde{\omega}$ is, indeed, an extreme point of $\mathbb{K} (\widetilde{\lambda}, \widetilde{\mu}).$ $\square$

\vskip0.2in
\noindent{\bf Example 3.4 \quad} Let $\lambda$ be the uniform distribution in the unit interval $[0,1].$ We shall use Lemma 3.3 and construct nongraphic extreme points of $\mathbb{K} (\lambda, \lambda)$ which are distributions in the unit square. To this end we start with the two points space $\mathbb{Z}_2 = \{0,1\}$ with the probability distribution $P$ where
$$P (\{0\}) = p, \,\,P (\{1\}) = q, \,\,0 < p < q < 1, \,\, p+q=1.$$
Now consider $\mathbb{Z}_2$-valued random variables $\xi, \eta$ with the joint distribution given by
$$P (\xi =0, \eta =0)=0, \,P (\xi=0, \eta=1) = P (\xi = 1, \eta = 0) = p,\, P(\xi=1, \eta = 1)=q-p. $$
Note that the joint distribution of $(\xi, \eta)$ is a nongraphic extreme point of $\mathbb{K} (P,P).$ Now consider an i.i.d sequence $\zeta_1, \zeta_2, \ldots$ of $\mathbb{Z}_2$-valued random variables independent of $(\xi, \eta)$ and having the same distribution $P.$ Put
$$\zta = (\zeta_1, \zeta_2, \ldots).$$
Then by Lemma 3.3 the joint distribution $\omega$ of $((\xi,\zta ),(\eta, \zta))$ is an extreme point of $\mathbb{K}(\nu, \nu)$ where $\nu = P \otimes P \otimes \ldots$ in $\mathbb{Z}_2^{\{0,1,2, \ldots\}}.$ Furthermore, since $(\xi, \eta)$ is nongraphic so is $((\xi,\zta ),(\eta, \zta)).$ Denote by $F_p$ the common probability distribution function of the random variables
$$\widetilde{\xi} = \frac{\xi}{2} + \sum_{j=1}^{\infty} \frac{\zeta_j}{2^{j+1}}, \quad \widetilde{\eta} = \frac{\eta}{2} + \sum_{j=1}^{\infty} \frac{\zeta_j}{2^{j+1}}. $$
Then $F_p$ is a strictly increasing and continuous function on the unit interval and therefore the correspondence $t \rightarrow F_p(t)$ is a homeomorphism of $[0,1].$ Put $\xi^{\prime} = F_p (\widetilde{\xi}), \eta^{\prime} = F_p (\widetilde{\eta}).$ Then the joint distribution $\omega$ of $(\xi^{\prime}, \eta^{\prime})$ is a nongraphic extreme point of $\mathbb{K}(\lambda, \lambda).$

Now we consider the case when $X_1$ and $X_2$ are finite sets, $G$ is a finite group acting on each $X_i,$ the number of $G$-orbits in $X_1, X_2$ and $X_1 \times X_2$ are respectively $m_1, m_2$ and $m_{12}$ and $\mu_i$ is a $G$-invariant probability distribution in $X_i$ with support $X_i$ for each $i=1,2.$ For any probability distribution $\lambda$ in any finite set denote by $S(\lambda)$ its support set. We first note that Theorem 3.2 assumes the following form.
\vskip0.2in
\noindent{\bf Theorem 3.5 \quad} A probability distribution $\omega \in \mathbb{K} (\mu_1, \mu_2)$ is an extreme point if and only if there is no nonzero real-valued function $\zeta$ on $S(\omega)$ satisfying the following conditions:
\begin{itemize}
 \item [(i)] $\zeta (gx_1, gx_2) = \zeta (x_1, x_2)$ $\qquad\forall$ $(x_1, x_2) \in S(\omega), g \in G;$
\item[(ii)] $\sum\limits_{x_{2} \in X_{2}} \zeta (x_1, x_2) \omega (x_1, x_2) = 0$ $\qquad\forall$ $x_1 \in X_1;$
\item[(iii)] $\sum\limits_{x_{1} \in X_{1}} \zeta (x_1, x_2) \omega (x_1, x_2) = 0$ $\qquad \forall$ $x_2 \in X_2.$ 
\end{itemize}
\vskip0.1in
\noindent{\bf Proof \quad} Immediate. $\square$
\vskip0.2in
\noindent{\bf Corollary 3.6 \quad} Let $\omega_1, \omega_2$ be extreme points of $\mathbb{K}(\mu_1, \mu_2)$ and $S(\omega_1) \subseteq S (\omega_2).$ Then $\omega_1 = \omega_2.$ In particular, any extreme point $\omega$ of $\mathbb{K}(\mu_1, \mu_2)$ is uniquely determined by its support set $S(\omega).$
\vskip0.2in
\noindent{\bf Proof \quad} Suppose $\omega_1 \neq \omega_2.$ Then put $\omega = \frac{1}{2} (\omega_1 + \omega_2).$ Then $\omega \in \mathbb{K} (\mu_1, \mu_2)$ and $\omega$ is not an extreme point. By Theroem 3.5 there exists a nonzero real-valued function $\zeta$ satisfying conditions (i)-(iii) of the theorem. By hypothesis $S(\omega) = S(\omega_2).$ Define
$$\zeta^{\prime} (x_1, x_2) = \frac{\zeta (x_1, x_2) \omega (x_1, x_2)}{\omega_2 (x_1, x_2)} \quad \mbox{where}\quad (x_1, x_2) \in S (\omega_2). $$
Then conditions (i)-(iii) of Theorem 3.5 are fulfilled when the pair $\zeta, \omega$ is replaced by $\zeta^{\prime}, \omega_2$ contradicting the extremality of $\omega_2.$ $\square$
\vskip0.2in
\noindent{\bf Corollary 3.7 \quad} For any $\omega \in \mathbb{K} (\mu_1, \mu_2)$ let $N(\omega)$ denote the number of $G$-orbits in its support set $S(\omega).$ If $\omega$ is an extreme point of $\mathbb{K}(\mu_1, \mu_2)$ then
$$\max (m_1, m_2) \leq N(\omega) \leq m_1 + m_2. $$
In particular, the number of extreme points in $\mathbb{K}(\mu_1, \mu_2)$ does not exceed
$$\sum_{\max (m_{1},m_{2}) \leq r \leq m_{1} + m_{2}} {m_{12} \choose r}. $$
\vskip0.2in
\noindent{\bf Proof \quad}  Let $\omega$ be an extreme point of $\mathbb{K}(\mu_1, \mu_2).$ Suppose $N(\omega) > m_1 + m_2.$ Observe that all $G$-invariant real-valued functions on $S(\omega)$ constitute a linear space of cardinality $N(\omega).$ Functions $\zeta$ satisfying conditions (i)-(iii) of the theorem constitute a subspace of dimension $\ge N(\omega) - (m_1 + m_2),$ contradicting the extremality of $\omega.$ For any distribution $\omega$ in $\mathbb{K}(\mu_1, \mu_2)$ we have $N(\omega) \ge m_i,$ $i=1,2.$ This proves the first part. The second part is now immediate from Corollary 3.6. $\square$ 

\vskip0.2in
\noindent{\bf Corollary 3.8\,\,(Birkhoff-von Neumann Theorem)\quad}  Let $X_1 = X_2 = X,$ $\# X = m,$ $\mu_1 = \mu_2 = \mu$ where $\mu (x) = \frac{1}{m}$ $\forall$ $x \in X.$ Then any extreme point $\omega$ in $\mathbb{K}(\mu, \mu)$ is of the form
$$\omega (x,y) = \frac{1}{m} \delta_{\sigma (x) y} \quad \forall \,\,\,x, y \in X $$
where $\sigma$ is a permutation of the elements of $X.$

 \vskip0.2in
\noindent{\bf Proof \quad}  Without loss of generality we assume that $X = =\{1,2,\ldots,m \}$ and view $\omega$ as a matrix of order $m$ with nonnegative entries with each row or column total being $1/m.$ First assume that in each row or column there are at least two nonzero entries. Then $\omega$ has at least $2m$ nonzero entries and by Corollary 3.7 it follows that every row or column has exactly two nonzero entries. We claim that for any $i \neq i^{\prime},$ $j \neq j^{\prime}$ in the set $\{1,2,\ldots,m\}$ at least one among $\omega_{ij}, \omega_{ij^{\prime}}, \omega_{i^{\prime}j}, \omega_{i^{\prime}j^{\prime}} $ vanishes. Suppose this is not true for some $i \neq i^{\prime},$ $j \neq j^{\prime}.$ Put
$$p = \min \left \{\omega_{rs} \big | (r,s): \omega_{rs} > 0 \right \}. $$
Define
$$\omega_{rs}^{\pm} = \left \{\begin{array}{l}
 \omega_{rs} \pm p \quad \mbox{if} \quad r=i, s=j \quad \mbox{or} \quad r=i^{\prime}, s=j^{\prime},\\
\omega_{rs} \mp p \quad \mbox{if} \quad r=i^{\prime}, s=j \quad \mbox{or} \quad r=i, s=j^{\prime},\\
\omega_{rs} \qquad \mbox{otherwise.}\end{array} \right . $$
Then $\omega^{\pm} \in \mathbb{K}(\mu, \mu),$ $\omega^{+} \neq \omega^{-}$ and $\omega = \frac{1}{2} (\omega^{+} + \omega^{-}),$ a contradiction to the extremality of $\omega.$ Now observe that permutation of columns as well as rows of $\omega$ lead to extreme points of $\mathbb{K}(\mu, \mu).$ By appropriate permutations of columns and rows $\omega$ reduces to a tridiagonal matrix of the form
$$\widetilde{\omega} = \left [\begin{array}{ccccccccc}
p_{11} & p_{12} & 0 & 0 &\ldots &\ldots&\ldots&\ldots& 0 \\ 
p_{21} & 0& p_{23}&0&\ldots&\ldots &\ldots&\ldots&0 \\  
0 &p_{32}&0&p_{34}&\ldots&\ldots &\ldots&\ldots& 0\\ 
\ldots&\ldots&\ldots&\ldots&\ldots&\ldots &\ldots&\ldots&\ldots \\ 
0&0&\ldots&\ldots&\ldots&0 &p_{n-1\,n-2}&0&p_{n-1\,n} \\ 
0&0&\ldots&\ldots&\ldots& 0&0&p_{n\,n-1}& p_{nn}  
\end{array}
 \right ] $$
where the $p$'s with suffixes are all greater than or equal to $p.$ Now consider the matrices
$$\lambda^{\pm} = \left [\begin{array}{cccccc}
p_{11} \pm p &p_{12} \mp p&0&0&0& \ldots \\
p_{21} \mp p&0&p_{23}\pm p &0&0& \ldots \\
0 & p_{32}\pm p &0& p_{34} \mp p&0&\ldots \\
\ldots &\ldots &\ldots &\ldots &\ldots &\ldots  
\end{array}
 \right ] $$
Then $\lambda^{\pm} \in \mathbb{K} (\mu, \mu)$ and $\widetilde{\omega} = \frac{1}{2} (\lambda^{+} + \lambda^{-}),$ contradicting the extremality of $\widetilde{\omega}$ and therefore of $\omega.$ In other words any extreme point $\omega$ of $\mathbb{K}(\mu, \mu)$ must have at least one row with exactly one nonzero entry. Then by permutations of rows and columns $\omega$ can be brought to the form
$$\omega_1 = \left [\begin{array}{c|cc} 
 1/m & 0 & 0 \ldots 0 \\ \hline
0 && \\
\vdots  && \widehat{\omega} \\
0 && \end{array} \right ] $$
where $\frac{m}{m-1} \widehat{\omega}$ is an extreme point of $\mathbb{K} (\widehat{\mu}, \widehat{\mu})$ where $\widehat{\mu}$ is the uniform distribution on a set of $m-1$ points. Now an inductive argument completes the proof. $\square$

We conclude with the remark that it is an interesting open problem to characterize the support sets of all extreme points of $\mathbb{K} (\mu_1, \mu_2)$ in terms of $\mu_1$ and $\mu_2.$


\newpage
\begin{thebibliography}{99}
\bibitem{rbbtesr}\label{rbbtesr} R. B. Bapat and T. E. S. Raghavan, {\it Nonnegative Matrices and Applications,} Cambridge University Press, Cambridge 1997. 

\bibitem{rbh}\label{rbh} R. Bhatia, {\it Matrix Analysis,} Springer Verlag, New York 1996.

\bibitem{gb}\label{gb} G. Birkhoff, {\it Tres observaciones sobre el algebra lineal,} University. Nac. Tucuman Rev, Ser. A5 (1946) 147-150.

\bibitem{krp}\label{krp} K. R. Parthasarathy, {\it Extreme points of the convex set of stochastic maps on a $C^{\ast}$ algebra,} Inf. Dim. Analysis, Quantum Probab. and Rel. Topics, 1 (1998) 599-609.

\bibitem{wfs}\label{wfs} W. F. Stinespring, {\it Positive functions on $C^{\ast}$ algebras,} Proc. Amer. Math. Soc., 6 (1955) 211-216. 

\bibitem{jvn}\label{jvn} J. von Neumann, {\it A certain zero-sum two-person game equivalent to an optimal assignment problem,} Ann. Math. Studies, 28 (1953) 5-12.
\end{thebibliography}
\end{document}